\documentclass{amsart}
\usepackage{amsmath,amssymb,amscd,mathrsfs, amsthm,multirow, graphicx}

\newtheorem{thm}{Theorem}
\newtheorem{prop}[thm]{Proposition}
\newtheorem{cor}[thm]{Corollary}

\theoremstyle{definition}
\newtheorem{define}[thm]{Definition}

\def\A{\mathbb A}

\def\P{\mathbb P}

\def\I{\mathscr I}

\def\O{\mathcal O}
\def\T{\mathcal T}
\def\NN{\mathcal N}
\def\X{\mathcal X}

\def\mm{\mathfrak m}

\def\vp{\varphi}

\def\ra{\rightarrow}

\def\>{\rangle}
\def\<{\langle}

\begin{document}
\title[Approximation for degree 2 del Pezzos]{Weak Approximation for General Degree Two Del Pezzo
Surfaces }
\author{Amanda Knecht}
%\titlerunning{Approximation for degree 2 del Pezzos}

\begin{abstract} We address weak approximation for
certain del Pezzo surfaces defined over the function field
of a curve. We study the rational connectivity of the smooth locus
of degree two del Pezzo surfaces with two  $A_1$ singularities in order to prove weak approximation for degree two del Pezzo surfaces with square-free discriminant.
\end{abstract}

\maketitle

\tableofcontents

\section{Introduction}
A standard question in Arithmetic Geometry is, ``Does a variety $X$
defined over a number field $K$ contain $K$-rational points, and if
so,  how are these points distributed on $X$?''   We say that $X$ satisfies weak
approximation if for any finite set of places of $K$ and points of
$X$ over the completion of $K$ at  these places, there exists a $K$-rational point of $X$
which is arbitrarily close to these points.  In this paper, we study
varieties defined over the function field of a smooth curve instead
of a number field.  In this context , rational points correspond to sections of fibrations over a curve, and proving weak approximation corresponds to finding sections with prescribed jet data in a finite number of fibers.

The existence of sections of rationally connected fibrations was proven by Graber, Harris and Starr \cite{GHS}.  Koll\'ar, Miyaoka, and Mori proved the existence of sections through a finite set of prescribed points in smooth fibers   \cite{K} IV.6.10 and \cite{KMM}  2.13.  The existence of sections with prescribed finite jet data through smooth fibers, i.e. weak approximation at places of good reduction, was proven by Hassett and Tschinkel \cite {WA}. In the same paper Hassett and Tschinkel conjectured that  a smooth rationally connected variety
$X$ defined over the function field of a smooth curve satisfies weak
approximation even at places of bad reduction.

There are a few cases where weak approximation over the function field of a curve is known to hold at all places.  Colliot-Th\'el\`ene and Gille proved weak approximation holds for stably rational varieties,
connected linear algebraic groups and homogeneous spaces for 
these groups,homogeneous space fibrations over varieties that satisfy weak approximation \cite{CTG}. In particular they prove that   conic bundles over $\P^1$ and del Pezzo surfaces of degree at least four satisfy weak approximation at all places.  The cases of del Pezzo surfaces of degree less that four are still open.  It is known that cubic surfaces with square-free discriminant satisfy weak approximation even at places of bad reduction \cite{WABR} .  This paper address weak approximation at places of bad reduction for degree 2 del Pezzo surfaces.

Let $k$ be an algebraically closed field of characteristic zero, $B$
a smooth curve over $k$ with function field $K = k(B)$, and $\overline{B}$ the smooth projective model of $K$ with $S:=\overline{B}/B$.

\begin{thm}Let $X$ be a smooth degree two del Pezzo surface over $K$, and $\pi: \X \ra B$ a proper model of $X$ (i.e. an algebraic space over flat over $B $ with generic fiber $X$.) Suppose the singular fibers of $\pi$ are degree two log del Pezzo surfaces with at most two $A_1$ singularities.  Let  $\X^{sm}$ be the locus where $\pi$ is smooth.  If there exists a section $s : B \ra \X^{sm}$, then the sections of $\X^{sm} \ra B$ satisfy approximation away from $S$.
\end{thm}

The proof of Theorem 1 relies on approximation results of \cite{WABR} and extending the rational connectivity results of \cite{KM} to show that in $X$ is a degree two log del Pezzo surface with one $A_1$ singularity, then each smooth point of $X$ is contained in a free proper rational curve contained in the smooth locus of $X$.   This proof also shows that  the sections of $\X^{sm} \ra B$ satisfy approximation away from $S$ when the singular fibers are degree two del Pezzo surfaces  with  the following singularity types:
$$ A_2, 2A_2,A_2+A_1,  A_3, A_3 + A_2 , A_3 + A_1, A_4 , A_4 + A_1, A_4 + A_2, A_5, A_6, D_4 .$$

The rationality of $X$ and properness of $\pi$ imply the existences of sections of $\pi$ \cite{GHS}.  When the model is regular, all sections of $\pi$ are contained in the smooth locus, so we find:

\begin{cor}
Let $X$ be a smooth degree two del Pezzo surface over $K$.  If $X$ admits a regular proper model $\pi: \X \ra B$ whose singular fibers are degree two log del Pezzo surfaces with two $A_1$ singularities, then weak approximation holds for $X$ away from $S$.  
\end{cor}
 This corollary is applicable whenever $B$ is a smooth projective curve because $X$ will always admit a regular proper model.  There exist smooth degree two del Pezzo surfaces that admit models with worse than $A_1$ singularities.  
 For example the family $$\{x^3 y - xy ^3 + tz^4 = w^2\} \subset \P(1,1,1,2) \times \A^1_t$$
over the $t$-line has a fiber with worse than even rational double points when $t=0$. But, Corollary 2 says weak approximation holds for `generic' degree two del Pezzo surfaces.

Let $X$ be a normal projective surface of degree 2 over $k$ with
ample anti-canonical sheaf $\omega_X^{-1}$.  Furthermore, suppose
that the singularities of $X$ are at most rational double points. We
call $X$ a degree 2 log del Pezzo surface with rational double
points. Such an $X$ is a double cover of $\P^2$ branched over a
quartic curve without multiple components thus may be thought of as a quartic hypersurface in the weighted projective space  $\P(1,1,1,2)$  defined by the equation $w^2=f_4(x,y,z)$ where $f_4 \in k[x,y,z]_4$ \cite{HW} .

Assume that $X$ is a smooth degree two del Pezzo surface defined over a smooth curve $B$, and $f:\X \ra B$ is a model of $X$ with smooth total space.  We would like to define the discriminant for this model, so we replace the singular space $\P(1,1,1,2)$ with its minimal resolution $W$.  Now the fibers of $f$ are degree four hypersurfaces in $W$.  Suppose $Y$ is a degree four hypersurface in $\P(1,1,1,2)$ and contains $[0,0,0,1]$.  Let $\widetilde{Y}$ be the proper transform of $Y$ in $W$.  Then $\widetilde{Y}$ will intersect the exceptional $\P^2$  along a conic of self-intersection $-4$.

Let $df: \T_{\X} \ra f^{\ast} \T_B$ be the differential map from a rank three vector bundle to a rank one vector bundle over $\X$.  Let $D(df)$ be the  degeneracy locus  of $df$  in $\X$ and define the discriminant $\Delta$  of $X$ to be the image of $D(df)$ in B.  The expected dimension of $D(df)$ is zero, but the cone point $[0,0,0,1]\in \P(1,1,1,2)$ can lead to components of dimension one as seen in the previous paragraph.  The excess intersection formula \cite{Ful} 6.3 tells us that the one dimensional components of $D(df)$ have multiplicity two in $\Delta$.  We also know that in the multiplicity of $\Delta$ at  each zero dimensional component of $D(df)$ is the sum of the Milnor numbers of the singularities in the corresponding fiber \cite{Ful} 7.1.14.  We say that the discriminant of $X$ is square-free if each component has multiplicity one, i.e. the  fibers of $f: \X  \ra B$ each have at most one $A_1$ singularity.
   
  \begin{cor}Suppose $X$ is a smooth degree two del Pezzo surface  defined over a smooth curve $B$.  If the discriminant of $X$ is square-free, then  sections of 
$f:\X \ra B$ satisfy approximation away from $S$.
\end{cor}

\noindent \textbf{Acknowledgments:} I am grateful to my thesis advisor Brendan Hassett for the many conversations we have had about this topic.  I also benefitted from talking to Damiano Testa and Brad Duesler. 
This material is based upon work supported by the National Science Foundation under Grants 0134259 and 0240058.

\section{Weak Approximation}

Let $K$ be a number field or the function field of a 
curve $B$. Let  $S$ be a finite set 
of places of $K$ that  contains the archimedean places .  For each place $\nu$ of $K$, let $K_{\nu}$ denote the $\nu$-adic completion of $K$. Let $X$ be an algebraic variety defined over $K$ and  $X (K )$ the set of $K$-rational 
points. One says that weak approximation holds for $X$ away from $S$  if for any finite set of 
places 
$\{\nu_i \}_{i \in I}$  of $K$, $\nu_i \notin S$ for $i \in I$, and $\nu_i$ -adic open subsets $U_i \subset X ( K_{\nu_i})$, there is a rational 
point $x\in X(K)$ such that its image in each $X ( K_{\nu_i})$  is contained in $U_i$ . 

By restricting ourselves to the case of function  fields, we can formulate a more geometric description of weak approximation.  Let $B$ be a smooth curve over an algebraically closed field of characteristic zero   with function field $K = k(B)$.  Let $\overline{B}$  be a smooth projective model of $K$ and set $S=\overline{B}/B$.  Let $X$ be a smooth proper variety over $K$, $\pi: \X \ra B$ a proper flat model (existence was proven in \cite{N}), and $X^{sm}$ the smooth locus of $\pi$.   Since $\pi$ is a proper morphism, sections $s: B \ra \X$ of $\pi$ correspond to $K$-rational points of $X$.  The analogue of local points are the $N$-jets defined below.

\begin{define} Let $\pi: \mathcal{X} \ra B$ be a proper model of
$X$ over $B$:
\begin{itemize}
\item[-] An \emph{admissible section of $\pi$} is a section
$s: B \ra \mathcal{X}^{sm}$.
\item[-] An \emph{admissible $N$-jet of $\pi$ at $b\in B$} is a section of
$$\mathcal{X}^{sm} \times_B \textrm{Spec}(\O_{B,b}/\mm_{B,b}^{N+1}) \ra
\textrm{Spec}(\O_{B,b}/ \mm_{B,b}^{N+1})$$ whose image is a smooth
point of $\X_b$.
\item[-] An \emph{approximable $N$-jet of $\pi$ at $b\in B$} is a
section of $$\mathcal{X} \times_B
\textrm{Spec}(\O_{B,b}/\mm_{B,b}^{N+1}) \ra \textrm{Spec}(\O_{B,b}/
\mm_{B,b}^{N+1})$$ that may be lifted to a section of $\hat{\X_b}
\ra \hat{B_b}$ where $\hat{B_b}:=$ Spec$(\hat{\O_{B,b}})$ and
$\hat{\X_b}:= \X \times_B \hat{B}_b$.
\end{itemize}
\end{define}
\noindent Note that Hensel's Lemma implies that every admissible
$N$-jet is approximable, and every section is admissible when the model is regular.  We can now formulate a geometric notion of weak approximation:

\begin{define} We say that \emph{$X$ satisfies weak approximation} away from $S$ if any finite collection of approximable jets of $\pi$ can be realized by a section $s: B \ra \X$. If $\X$ is a regular model this is equivalent to the condition that any 
collection of admissible jets of $\pi$ can be realized by a section 
$s : B \ra X^ {sm}.$\\
When we want to refer to a specific model $\pi: \X \ra B$, we  say that \emph{sections of $\pi$  satisfy approximation} away from $S$ if any finite collection of approximable jets of $\pi$ can be realized by a section $s: B \ra \X$.
\end{define}

\section{Notions of Rational Connectivity}
\begin{define}  A variety $X$ is
\emph{rationally connected} if there is a family of proper algebraic  curves $g:U \ra Y$ whose fibers are irreducible rational curves  and a cycle morphism $u: U \ra X$ such that 
$$u^{2}:\;\;\;     U \times_Y U  \longrightarrow X \times X $$
is dominant.
\end{define}
When $X$ is defined over an uncountable algebraically closed field $k$, rational connectivity is equivalent to the condition that any two very general  points  in $X(k)$ can be joined by an irreducible projective rational  curve contained in $X$ \cite{K} IV.3.6. 
\begin{define} 
Let $X$ be a smooth algebraic variety and
$f : \P^1 \ra X$ a nonconstant morphism, so we have an isomorphism
$$f^{*}\T_X \simeq \O_{\P^1}(a_1)\oplus \cdots \oplus  \O_{\P^1}(a_{dim X})$$ for suitable
 integers $a_1,\ldots , a_{dim X}.$
  Then $f$ is \textit{free} (resp. \textit{very free}) if each $a_i \geq 0$ (resp. $a_i \geq 1$).
\end{define}

When $X$ is a smooth variety over an algebraically closed field, being rationally connected is equivalent to containing one very free curve \cite{K} IV 3.7.  Suppose we are in this situation.  Then there exists a unique largest nonempty subset $X^0 \subset X$ such that for any finite collection of distinct points in $X$ there is a very free rational curve contained in $X^0$ which contains these points as smooth points. Moreover, any rational curve that meets $X^0$ is contained in $X^0$ \cite{K} IV.3.9.4.  There are no known examples where $X \neq X^0$.  The case of $X=X^0$ leads to the following definition.

\begin{define}[\cite{WABR} 14] A smooth rationally connected
variety $X$ is \textit{strongly rationally connected} if any of the
following conditions hold:
\begin{enumerate}
  \item [(1)] for each point $x \in X$, there exists
a rational curve $f : \P^1 \ra X$ joining $x$ and a generic point in
$X$;
\item[(2)] for each point $x\in X$, there exists a very free rational curve containing $x$.
  \item [(3)]for any finite collection of points $x_1, \ldots, x_m \in X$, there
exists a very free rational curve containing the $x_j$ as smooth
points;
  \item [(4)] for any finite collection of jets $$j_{N,i}\in \textrm{Spec}\frac{ k[\varepsilon]}{
   \left<\varepsilon ^{N+1}\right> }\subset X , \; \; i= 1, \ldots, m$$
supported at distinct points
  $x_1, \ldots,  x_m$, there
exists a very free rational curve smooth at $x_1,\ldots,x_m$ and
containing the prescribed jets.
\end{enumerate}\end{define}

We note here that in Definition 14 of \cite{WABR} property $(2)$ reads `free' and not `very free,'  but having a free curve through every point has not been proven to be equivalent to any of the other conditions.  When the variety $X$ is a surface, the word `very free' can be replaced by `free' in property $(2)$.
  \begin{prop}
A smooth rationally connected surface $X$ is strongly
rationally connected if and only if  for each point $x \in X$,
there exists a free rational curve containing $x$.
\end{prop}
\begin{proof} When there is a family of rational curves through $x$ dominating $X$,
one of these curves is guaranteed to be free \cite[1.1]{KMM}.    Now
suppose there is a free curve through every point in $X$. For any
$x\in X$, let $R$ be a free rational curve containing $x$.  Let
$D$ be a free curve containing the generic point. Since $R$ and $D$
are free, we know that $R^2, D^2 \geq 0$. By the Hodge Index Theorem
$D^2R^2- (D.R)^2 < 0$.  Thus the curves $D$ and $R$ intersect in
$X$.  We can smooth this curve to a rational curve
 joining $x$ and the generic point \cite[II.7.6]{K}.
\end{proof}

The main point of this paper is to prove that the smooth locus of a degree two log del Pezzo surface with an $A_1$ singularity is strongly rationally connected because Hassett and Tschinkel proved the following theorem:

\begin{thm} \textsc{ (\cite{WABR} 15) } Let $\pi : \X \ra  B$ be a smooth morphism whose fibers are 
strongly rationally connected. Assume that $\pi$ has a section. Then sections 
of $\X \ra B$ satisfy approximation away from $S$ .
\end{thm}

\section{Combs and Free Curves}
Our main goal of this paper is to find a free rational  curve through any smooth point in $X$.   Suppose we have a surface containing many very free rational curves.  Suppose also that we have found a  rational curve through a point $x\in X$ which is not free.  This section shows how we use combs to create a free rational curve through the point $x$.

\begin{define}A \textit{comb with $p$ teeth} is a
connected and reduced  curve with $p + 1$ irreducible components $D,C_1, \ldots
, C_p$  such that:\begin{enumerate}
\item $D$ is smooth and called the \textit{handle} ;
\item every $C_i$ is isomorphic to $\P^1$ and these are calles the \textit{teeth};
\item the only singularities of $C$ are ordinary nodes;
\item every $C_i$ intersects $D$ in a single point, and $C_i\cap C_j=\emptyset $ when $i\neq j$.
\end{enumerate}
 A \textit{rational comb} is a comb whose
handle is a smooth rational curve.
\begin{figure}[h!]\centering
\includegraphics[height=2in]{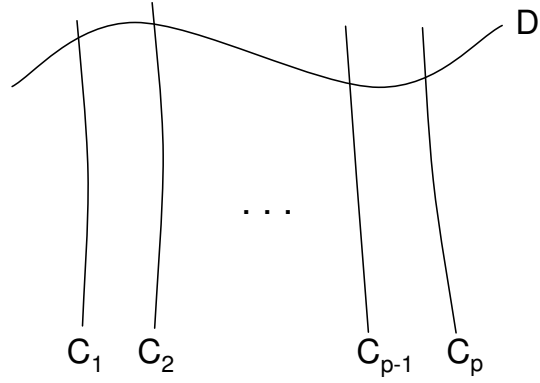}
\caption{Comb with p teeth}\label{fig:Comb}
\end{figure}
\end{define}

We construct a comb whose handle is an arbitrary nonfree rational
curve $D$ on a smooth quasi-projective variety $X$.  Assume
there exists several very free rational curves $C_1, . . . ,C_p \subset
X$ meeting $D$ at distinct smooth points $d_i \in D$. We want  to deform the union $D \cup  C_1 \cup \ldots \cup C_p$ into a
very free rational curve.  One construction of such a comb was considered
by Koll\'{a}r, Miyaoka, and Mori, see \cite[II.7.9, II.7.10]{K}. A
more recent construction of a comb that deforms nicely comes from
the techniques of Graber, Harris, and Starr.  This construction
considers deformations in the Hilbert scheme of $X$.  We include
here some essential facts about Hilbert schemes which are used  in
constructing the comb that frees $D$.

There is a $k$-scheme Hilb$(X)$ parameterizing closed subschemes of
$X$.  Assume  $X$ is smooth, and let $Z$ be a closed subscheme of
$X$ given by the ideal sheaf $\I_Z$. Assume  $Z$ is a local complete
intersection, i.e.  $\I_Z$ is locally generated by codim$(Z,X)$
elements. Then $\I_{Z}/\I_{Z}^{2}$ is  locally free
 on $Z$ and its dual is the normal bundle $\NN_{Z/X}$ of $Z$
in $X$. In this situation the following hold \cite[I.2.8]{K}:
\begin{itemize}
\item[1.] $\T_{[Z]} $Hilb$(X) \cong $H$^0(Z,\NN_{Z/X})$.
\item[2.] If H$^1(Z,\NN_{Z/X}) = 0$, then Hilb$(X)$ is smooth at
$[Z]$.
\end{itemize}

Suppose there are disjoint smooth rational curves $C_1,\ldots ,C_p$
each meeting a rational curve $D$ transversely in a single point
$d_i$, $i=1 , \ldots, p$. Fix a point $P\in D$ that is not one of
the nodes $d_i$.  Then we can deform their union $D^{\ast} = D \cup
C_1\cup \ldots \cup C_p$ into an irreducible curve provided the
following conditions hold \cite[2.6]{GHS}:
\begin{itemize}
\item[1.] The sheaf $\NN_{D^{\ast}/X}$ is generated by global sections.
\item[2.] H$^1(D^{\ast},\NN_{D^{\ast}/X}) = 0$.
\end{itemize}
We can see this by considering the sequence $$ 0 \ra \NN_{D/X}
 \ra {\NN_{D^{\ast}/X}}|_{D} \ra Q \ra 0$$
where $Q$ is a torsion sheaf supported at the points $d_i = D \cap
C_i$. If $\NN_{D^{\ast}/X}$ is generated by global sections, then we
can find a global section $s \in$ H$^0(C,\NN_{C/X})$ such that, for
each $i=1, \ldots p$, the restriction of $s$ to a neighborhood of
$d_i$ is not in the image of $\NN_{D/X}$. This means  $s$
corresponds to a first-order deformation of $C$ that smoothes the
nodes $d_i$ of $D^{\ast}$. If the second condition holds, there is
no obstruction to lifting our first-order deformation to a global
deformation of $D^{\ast}$ smoothing the nodes.

For a general fixed  point $P\in D$, we would like to create a comb
which deforms into an irreducible free curve while still passing
through $P$.  In order to achieve this, we replace the conditions
above with  new conditions:
\begin{itemize}
\item[1'.] The sheaf $\NN_{D^{\ast}/X}(-P)$ is generated by global sections.
\item[2'.] H$^1(C,\NN_{D^{\ast}/X}(-P)) = 0$.
\end{itemize}

\begin{thm} \textsc{ (\cite{KA} 27) } \label{thm:KA}
 Let $X$ be a smooth projective variety of dimension at least 3
over an algebraically closed field. Let $D\subset X$ be a smooth
irreducible curve and $\mathscr{M}$ a line bundle on $D$. Let $C
\subset X$ be a very free rational curve intersecting $D$ and let
$\mathcal{C}$ be a family of nodal rational curves on $X$
parameterized by a neighborhood of $[C]$ in Mor$(\P^1,X)$. Then
there are curves $C_1, \ldots ,C_p \in \mathcal{C}$ such that
$D^{\ast} = D \cup C_1 \cup \ldots \cup C_p$ is an immersed  comb
and satisfies the following conditions:
\begin{itemize}
\item[1.] The sheaf $\NN_{D^{\ast}/X}$ is generated by global sections.
\item[2.] H$^1(D^{\ast},\NN_{D^{\ast}/X}(\mathscr{M}^{\ast})) = 0$, where $\mathscr{M}^{\ast}$ is the unique line bundle on
$D^{\ast}$ that extends $\mathscr{M}$ and has degree 0 on the $C_i$.
\end{itemize}
\end{thm} 
Thus on a variety of dimension at least three containing many very
free curves, rational curves meeting the very free curves can be
freed. We can generalize this result a little.  

\begin{prop}  Let $X$ be a smooth surface over an
algebraically closed field.  Let $D \subset X$ be an irreducible
rational curve and $P$ a point on $D$.  Let $C \subset X$ be a very free
rational curve intersecting $D$ and let $\mathcal{C}$ be a family of
rational curves on $X$ parameterized by a neighborhood of $[C]\in
$Mor$(\P^1,X)$. Then there are curves $C_1, \ldots ,C_p \in
\mathcal{C}$ such that $D^{\ast} = D \cup C_1 \cup \ldots \cup C_p$
is a comb whose nodes can be smoothed to create a free rational
curve containing $P$.
\end{prop}

\begin{proof}   We  modify the smoothness of $D$  and dimension of $X$ in Theorem \ref{thm:KA}.\\
 \emph{ Step 1.}  Reduction to the case where
$D$ is a smooth rational curve.

 Suppose  $D$ is singular. Then we
will work with the normalization $\widetilde{D}\ra D$.  When
considering the Hilbert Scheme, it is essential that $D$ be embedded
in $X$, so we take an embedding $\widetilde{D} \hookrightarrow \P^3$
and consider the diagonal map $\Delta : \widetilde{D} \ra X \times
\P^3.$ This is an embedding with image $D'$ isomorphic to
$\widetilde{D}$.  Let $\pi_1$ denote the projection from $  X \times
\P^3 $ to the first factor.  Then $\pi_1(D')\simeq D$. Any moving of
$D'\subset  X \times \P^3$ can be projected to $X$ to give a
deformation of $D$. From now on assume $D\subset X$ is a smooth,
irreducible, rational curve.

\noindent\emph{ Step 2.}  Reduction from  surfaces to
three-folds.

 Suppose $S$ is a smooth surface and $f: \P^1 \ra S$
is a morphism.  Consider the three-fold $X=\P^1 \times S$ and the
map $g=(id, f): \P^1 \ra X$. If we let $\pi_2: X \ra S$ be the
projection onto the second factor, then $f=\pi_2 \circ g.$  Since
$\pi_2$ is smooth, the map $g^{\ast} \T_X \ra f^{\ast}\T_S$ is
surjective . Thus if $g^{\ast} \T_X$ is
ample, then $f^{\ast}\T_S$ is ample.

\end{proof}

Now once we find a rational curve in $X^{sm}$ containing a given
smooth point $x$, we know how to free it.  So what is left  to do is
to find rational curves through points.

\section{Proof of the Main Theorem}
Keel and McKernan proved the smooth locus  of any  log del Pezzo surface is rationally connected \cite{KM}.  As stated in the Section 3, the main contribution made by this paper is the following refinement for degree two del Pezzos:  

\begin{thm} Let $X$ be a  degree two log del Pezzo surface with two $A_1$ singularities  and $X^{sm}$ the smooth locus of $X$.  Then $X^{sm}$ is strongly rationally connected.
\end{thm}

 \begin{proof}
 Let $ \vp: X \ra \P^2$ be a double cover of $\P^2$ branched over a quartic plane curve $Q$.  Since $X$ has two $A_1$ singularities, $Q$ has two nodes $q_1, q_2$ \cite{Barth} III.7.1. 
In order to prove the theorem, for every point $x\in X^{sm}$ we will find a free curve $C \subset X^{sm}$ containing $x$.

Consider arbitrary $x\in X^{sm}$.  Let $\vp(x)=p \in \P^2 \setminus \{q_1, q_2\}$. We must find a proper curve in $p \in \P^2 \setminus \{q_1, q_2\}$ which lifts to a free rational curve in $X^{sm}$.  The simplest plane curves are lines, so we can try to find lines through $p$ that lift to rational curves.  By the Hurwitz formula 
 we see that a generic line in $\P^2$ will lift to 
a curve of genus one in $X$.  However, lines simply  tangent to $Q$ will lift to rational curves.  Thus, we consider lines tangent to the branch curve $ Q$.

First we consider the case when $p$ is not contained in $Q$.   Let $x\in X$ such that $\vp(x)=p\notin Q$.
Let $\widetilde{Q}$ be the normalization of $Q$ and consider the projection from $\pi_p: \P^2 \setminus \{p\} \ra
 \P^1$.  Then $\pi_p$ restricted to the the genus two curve $\widetilde{Q}$ is a
 four-to-one cover of the projective line with ramification index equal to
ten.  Since the ramification number is ten,  there are at least four lines through $p$ that are tangent to $Q$ and three of these lines will not contain the singular point $q$.  We only need to pick one of these lines and analyze the four ways a line can be tangent to a quartic:  one simple tangent; one three tangent;  one four-tangent; or two simple tangents.  

Next  suppose that $p$ is a smooth point on the quartic $Q$.  If the tangent line of $Q$ at $p$ does not contain the singular point, $q$, then we can use the tangent line at $p$ to find a free curve through $x$ with the analysis below.  The tangent line at $p$ passing through through $q$ is not actually a problem.  If we consider projection from  $p$ again, we have a three to one cover of the line with ramification index equal to eight and are guaranteed that $p$
 is contained in at least four lines tangent to $Q$.  We simply  pick one of the lines missing the singular point.  
 
Now we analyze the four cases of tangency between a line an a quartic cure, and construct free curves in $X$ from lifts of the lines. \\

\noindent \textbf{One Simple Tangent.}
\begin{figure}[h!]\centering
\includegraphics[height=3in]{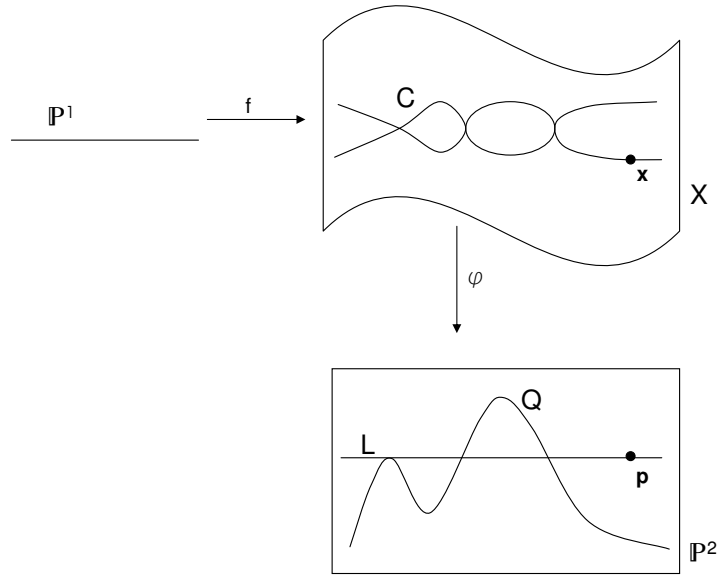}
\caption{One Simple Tangent}\label{fig:One Simple Tangent}
\end{figure}
 Suppose that we can find a line $L$ through $p$ that is
simply tangent to $Q$ at a smooth point of $Q$ and intersects the
quartic in two other distinct points. By the Hurwitz formula, we
know that the normalization, $\widetilde{C}$, of the curve $C
\subset X$ lying over $L$ is rational, thus $C$ is rational. Let $C$
be the image of the projective line under the map $f: \P^1 \ra X$.
\\We have an exact sequence
$$0 \longrightarrow \T_{\P^1} \longrightarrow f^{*}\T_X \longrightarrow
\NN_{f} \longrightarrow 0.$$  We know that $\T_{\P^1}=
\O_{\P^1}(2)$, and by the double point formula deg$\NN_f =C^2-
2(\textrm{\# nodes on } C)= 2-2=0$.  The extension $$0
\longrightarrow \O_{\P^1}(2) \longrightarrow f^{*}\T_X
\longrightarrow \O_{\P^1} \longrightarrow 0$$
 is given by an
element of Ext$^1(\O_{\P^1}, \O_{\P^1}(2))=$ H$^1(\P^1,
\O_{\P^1}(2))=0$. Thus, the sequence splits and  $f^{*}\T_X=
\O_{\P^1}(2)\oplus \O_{\P^1}(0)$, and
 $C$ is a free curve through the point $x$.

We can even create combs with $C$ as the handle in order to find very free curves through $x$.  Since $X^{sm}$ is rationally connected, we know that $X^{sm}$ contains a very free rational curve $D$. The Hodge Index Theorem  and the the fact the $D^2$ and $C^2$ are positive imply that  $D$ and $C$ intersect.  By Proposition 10, we see that there exists a very free curve through $x$. \\

\noindent\textbf{One Inflectional 3-Tangent.}
 Now assume  we can find a line $L$ through
$p$ that is three-tangent to $Q$ at a smooth point. As before, we let $C \subset X$ be the rational
curve lying over $L$ and containing $x$. With the same set up as the
previous case, we return to the exact sequence,
$$0 \longrightarrow \T_{\P^1} \longrightarrow f^{*}\T_X \longrightarrow
\NN_{f} \longrightarrow 0.$$  Again we find that $f^{*}\T_X=
\O_{\P^1}(2)\oplus \O_{\P^1}(0)$, and we have a free rational curve
through the point $x\in X$.\\\\

\noindent \textbf{One Inflectional 4-Tangent.}   Suppose we have a line $L$ through $p$   intersecting the quartic $Q$
at only one smooth point. Then
$L$ is four tangent to $Q$ at the point of intersection, and the
corresponding curve $C \subset X$ above $L$ will have a tacnode
above the intersection.  Because the arithmetic genus of the
normalization of the curve $C$ is -1, we know  $C$ must be the union
of two irreducible rational curves $C_1, C_2$. Instead of
considering the entire curve $C$, we shall focus on the component
$C_1$ passing through the original point $x \in X$.
\begin{figure}[h!]\centering
\includegraphics[height=1.7in]{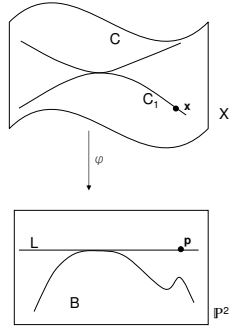}
\caption{One Inflectional 4-Tangent}\label{fig:4tangent}
\end{figure}

\noindent   Because each $C_i$ is mapped one-to-one onto the line
$L$ in $\P^2$ by the anti-canonical sheaf of $X$, we know
$C_i\centerdot K_X= -1$. Adjunction then implies that the self
intersection, hence the degree of the normal sheaf, of each $C_i$ is
-1. Hence, our curve $C_1$ is not free.

We use the dual curve of $Q$ to show that we can use a comb as in Proposition 12 to free $C_1$.  For each point $q\in Q$
let $\T_q Q$ be the tangent line to $Q$ at $q$.  Consider $\T_q Q$
as a point in the dual projective plane $(\P^2)^{\ast}$.  There is a
map from $Q$ to its dual curve $Q^{\ast}\subset (\P^2)^{\ast}$ which
sends a point $q$ to $\T_q Q$.  The singularities of $Q^{\ast}$ are:
nodes corresponding to bitangents to $Q$, cusps corresponding to
inflectional 3-tangents of $Q$, tacnodes corresponding to
inflectional 4-tangents of $Q$. The dual curve to a quartic with one node is a
planar curve of degree 10 and has only a finite number of
singularities. Thus, for an infinite number of points on the line
$L$, there are lines through these points that are simply tangent to
$Q$ and intersect $L$ transversely. These lines lift to free curves
in $X$ as shown above, and can be made into the very free curves that we need in order to apply the proposition.\\

\noindent \textbf{Two Simple Tangents.}
 \noindent  Finally suppose there exists a line $L$ through $p$ that is tangent
to $Q$ at two distinct nonsingular points.  Then the curve $C\subset
X$ lying above $L$ will have a node above each of these tangencies.
The arithmetic genus of the normalization of $C$ is -1, so we know
 $C$ is the union of two rational curves $C_1, C_2$.  As in the
case above, we find that the degree of the normal sheaf of each
$C_i$ is -1. Suppose $C_1$ is the component that contains the point
$x$. We can use very free curves intersecting $C_1$ to create a comb to
free $C_1$. The resulting curve will be a free rational curve
through $x$ as desired.
 \end{proof}

 \begin{figure}[h!]\centering
\includegraphics[height=1.7in]{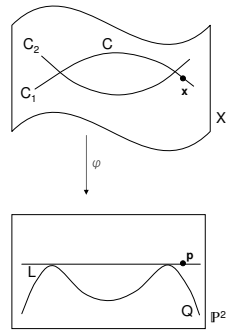}
\caption{Bitangent}\label{fig:Bitangent}
\end{figure}

 As stated in the introduction, this proof can be easily modified to show that the smooth locus of a degree two del Pezzo surface is strongly rationally connected when the singularity type of the surface is one of the following: 
$$ A_2, 2A_2,A_2+A_1,  A_3, A_3 + A_2 , A_3 + A_1, A_4 , A_4 + A_1, A_4 + A_2, A_5, A_6, D_4 .$$

\bibliography{KnechtWAdegree2}
	\bibliographystyle{plain}

\end{document}